\documentclass[11pt,twoside]{article}
\usepackage{amssymb}
\usepackage{graphicx}

\usepackage{graphicx,type1cm,eso-pic,color}
\usepackage{amssymb}
\usepackage{amssymb,amsmath}
\usepackage{graphicx,epsfig}
\usepackage{enumerate}
\usepackage{color}
\usepackage {multicol}

\numberwithin{equation}{section}
\newtheorem{theorem}{Theorem}
\newtheorem{lemma}{Lemma}

\newtheorem{corollary}{Corollary}
\newtheorem{example}{Example}

\makeatletter
\AddToShipoutPicture{
	\setlength{\@tempdimb}{.5\paperwidth}
	\setlength{\@tempdimc}{.5\paperheight}
	\setlength{\unitlength}{1pt}
	\put(\strip@pt\@tempdimb,\strip@pt\@tempdimc)
}
\makeatother

\usepackage{scalerel,stackengine}
\stackMath
\newcommand\reallywidehat[1]{%
	\savestack{\tmpbox}{\stretchto{%
			\scaleto{%
				\scalerel*[\widthof{\ensuremath{#1}}]{\kern-.6pt\bigwedge\kern-.6pt}%
				{\rule[-\textheight/2]{1ex}{\textheight}}
			}{\textheight}%
		}{0.5ex}}%
	\stackon[1pt]{#1}{\tmpbox}%
}
\begin{document}
	\setcounter{page}{1}

	\thispagestyle{empty}
	\markboth{}{}

	\pagestyle{myheadings}
	\markboth{ S.K.Chaudhary, N. Gupta }{ S.K.Chaudhary, N. Gupta }
	
	\date{}
	
	
	\noindent  
	
	\vspace{.1in}
	
	{\baselineskip 20truept
		
		\begin{center}
			{\Large {\bf General weighted cumulative residual (past) extropy of minimum (maximum) ranked set sampling with unequal samples}} \footnote{\noindent	{\bf * } Corresponding author E-mail: skchaudhary1994@kgpian.iitkgp.ac.in\\
				{\bf **} E-mail: nitin.gupta@maths.iitkgp.ac.in
			}			
		\end{center}

		\vspace{.1in}
		
		\begin{center}
			{\large {\bf Santosh Kumar Chaudhary*, Nitin Gupta**}}\\
			{\large {\it Department of Mathematics, Indian Institute of Technology Kharagpur, West Bengal 721302, India }}
			\\
		\end{center}

		\vspace{.1in}
		\baselineskip 12truept

		\begin{abstract}
			\noindent The general weighted cumulative residual extropy (GWCRJ) and general weighted cumulative past extropy (GWCPJ) are introduced in this paper. There are some results in relation to GWCPJ and GWCRJ. We take into account GWCRJ-based uncertainty measures for the minimal ranked set sampling technique with unequal samples (minRSSU). Additionally, we take into account GWCPJ-based uncertainty measures for the maximum ranked set sampling technique with unequal samples (maxRSSU). Stochastic comparison for Simple random sampling (SRS)  is discussed. We looked at the monotone properties of minRSSU and maxRSSU as well as stochastic comparison.
			Finally, two empirical estimators of GWCPJ and GWCRJ are obtained.
			
			\vspace{.1in}
			
			
		\end{abstract}
		
		\section{Introduction}\label{section1}
		The method of enhancing the accuracy of the sample mean as an approximation of the population means, McIntyre (1952) proposed ranked set sampling (RSS). The RSS process involves selecting $n$ sets of units, each of size $n$, at random from a population. The units in each set are then ranked visually or using some low-cost techniques. Pick the unit with the lowest rank out of the first $n$ units. Select the unit with the second-lowest rank from the second group of $n$ units. Up until selecting the unit with the highest rank in the $n$-th set, the process is repeated. Then, we measure the relevant variable using the units we have chosen. In many practical aspects related to applications in industrial statistics, environmental and ecological studies, statistical genetics, etc., RSS schemes have been playing a significant role. We assume that $ X_{SRS}=\{X_i, i=1,\dots, n\}$ denotes a simple random sampling (SRS) of size $n$ from a continuous distribution with probability density function (pdf) $f$ and cumulative
		distribution function (cdf) $F.$ In comparison to SRS, the RSS is a more effective sampling method for estimating the population mean. As a helpful improvement to the RSS technique, Qiu and Eftekharian (2021) investigated the extropy of the minimum and maximum ranked set sampling procedure with unequal samples (minRSSU and maxRSSU). In the minRSSU and maxRSSU, we draw $n$ simple random samples, in which the size of the $i-$th samples is $i, \  i = 1,\dots n.$ The one-cycle minRSSU or maxRSSU involves an initial ranking of $n$ samples of size $n$.  We collect $Z_i=X_{(1:i)i}$ for all $i=1,2,\dots, n, \ $ and $Y_i=X_{(i:i)i}$ for all $i=1,2,\dots, n, \ $ where $X_{(i:j)j}$ denotes the $i$-th order statistic from the $j$-th SRS of size $j.$ The resulting sample is called one-cycle minRSSU and maxRSSU of size $n$ and denoted by $Z_{minRSSU} =\{Z_i, i=1,\dots,n\}$ and $Y_{maxRSSU} =\{Y_i, i=1,\dots,n\},$ respectively (see Biradar and Santosha (2014)).

		Measures of information for RSS and its variations have been developed by several authors (see Eskandarzadeh et al.(2018), Raqab and Qiu (2019), Qiu and Raqab (2022) and Gupta and Chaudhary (2022) etc.). Qiu and Eftekharian (2021) considered the information content of minRSSU and maxRSSU in terms of extropy and also compared the extropy of these two sampling data with that of SRS and RSS data. Kazemi et al. (2021) considered uncertainty measures of minimum ranked set sampling procedure with unequal samples (MinRSSU) in terms of cumulative residual extropy and its dynamic version and they also compared the uncertainty and information content of cumulative residual extropy based on MinRSSU and SRS designs. Chacko and George (2022) considered the extropy properties of the RSS when ranking is not perfect and also obtained upper and lower bounds of the extropy of RSS.

		Extropy was introduced by Lad et al.(2015) as a measure of uncertainty as
		\begin{equation}
			J(X)=-\frac{1}{2} \int_{-\infty}^{\infty} f^2(x)dx .
		\end{equation}
		Jahanshahi et al. (2019) defined Cumulative residual extropy as
		\begin{equation}
			\xi (X) =-\frac{1}{2} \int_{-\infty}^{\infty} \bar{F}^2(x)dx .
		\end{equation}
		Cumulative past extropy is defined as 
		\begin{equation}
			\bar \xi (X) =-\frac{1}{2} \int_{-\infty}^{\infty} F^2(x)dx .
		\end{equation}
		Extropy and its various form have been studied by several authors (see Bansal and Gupta (2021), Gupta and Chaudhary (2022), and Kazemi et al.(2022) etc). Bansal and Gupta (2021) studied weighted extropies and past extropy of order statistics and k-record values. Gupta and Chaudhary (2022) studied the general weighted extropy of RSS where they defined general weighted extropy (GWJ) with weight $w(x)\geq 0$ as
		\begin{align*}
			J^{w}(X)&=-\frac{1}{2} \int_{-\infty}^{\infty}w(x) f^2(x)dx\\
			&=-\frac{1}{2}\mathbb{E}(\delta_X^{w}(U)),
		\end{align*}
		where $\delta_X^{w}(u)=w(F^{-1}(u))f(F^{-1}(u))$ and $U$ is uniformly distributed random variable on $(0,1)$, i.e., $U\sim$ Uniform$(0,1)$. This motivated us to generalize cumulative residual and past extropy and study it for ranked set sampling.

		Here, we defined general weighted cumulative past extropy (GWCPJ) of absolutely continuous random variable $X$ with weight $w(x)\geq 0$ as	
		\begin{align}\label{GWCPJdef}
			\bar{\xi}^w(X)&=-\frac{1}{2} \int_{-\infty}^{\infty} w(x) F^2(x)dx \nonumber\\
			&= -\frac{1}{2} \int_{0}^{1} w(F^{-1}(u)) \frac{u^2}{f(F^{-1}(u))}du \nonumber\\
			&=-\frac{1}{2}\mathbb{E}(\Lambda_X^{w}(U)),
		\end{align}
		where $\Lambda_X^{w}(u)= \frac{u^2 w(F^{-1}(u))}{f(F^{-1}(u))}.$
		
		\noindent Also, we define general weighted cumulative residual extropy (GWCRJ) of absolutely continuous random variable $X$ with weight $w(x)\geq 0$ as	
		\begin{align}\label{GWCRJdef}
			{\xi}^w(X)&=-\frac{1}{2} \int_{-\infty}^{\infty} w(x) \bar{F}^2(x)dx \nonumber\\
			&= -\frac{1}{2} \int_{0}^{1} w(F^{-1}(u)) \frac{(1-u)^2}{f(F^{-1}(u))}du \nonumber\\
			&=-\frac{1}{2}\mathbb{E}(\Delta_X^{w}(U)), 
		\end{align}
		where $\Delta_X^{w}(u)= \frac{(1-u)^2 w(F^{-1}(u))}{f(F^{-1}(u))}.$
		
		This paper is organized as follows:  section \ref{section2} discusses results on GWCPJ and GWCRJ. Section \ref{section3} provides examples and results for SRS of GWCPJ and GWCRJ. MinRSSU and maxRSSU for GWCRJ and GWCPJ, respectively are studied in section \ref{section4}. Stochastic comparison for MinRSSU and maxRSSU for GWCRJ and GWCPJ, respectively are obtained in section \ref{section5}. Section \ref{section6} is the monotone properties of MinRSSU and maxRSSU for GWCRJ and GWCPJ, respectively. In section \ref{section7}, we obtained empirical estimators of GWCRJ and GWCPJ. Finally, section \ref{section8} concludes this paper.
		
		\section{ Results on GWCPJ and GWCRJ}\label{section2}
		In the following theorem, we will study the behaviour of GWCPJ with respect to weight functions $w_1$ and $w_2$ and the dispersive order of random variables $X$ and $Y$.
		
		\begin{theorem}
			Let $X$  and $Y$ be nonnegative random variables with pdf's $f$ and $g$, cdf's $F$ and $G$, respectively having $u_X=u_Y<\infty$.\\
			\\
			(a) If $w_1$ is decreasing, $w_1(x)\leq w_2(x)$ and $X\le_{disp} Y$, then $\bar{\xi}^{w_1}(X)\ge \bar{\xi}^{w_2}(Y)$.\\
			(b)   If $w_1$ is decreasing, $w_1(x)\geq w_2(x)$ and $X\ge_{disp} Y$, then $\bar{\xi}^{w_1}(X)\le \bar{\xi}^{w_2}(Y)$.
		\end{theorem}
		\noindent \textbf{Proof}
		(a) Since $X\le_{disp} Y$, therefore we have $f(F^{-1}(u))\ge g(G^{-1}(u))$ for $u \in (0,1)$. Then using Theorem 3.B.13(b) of Shaked and Shanthikumar (2007), $X\le_{disp} Y$ implies that $X\ge_{st} Y$. Hence $F^{-1}(u) \ge G^{-1}(u)$ for all $u\in (0,1)$. Since $w_1$ is decreasing and $w_1(x)\leq w_2(x)$ , then $w_1(F^{-1}(u)) \le w_1(G^{-1}(u))\le w_2(G^{-1}(u))$. 
		Hence, 
		\begin{align}
			\Lambda_X^{w_1} (u)&=u^2 \frac{ w_1(F^{-1}(u))}{f(F^{-1}(u))}\nonumber\\
			&\le u^2 \frac{ w_2(G^{-1}(u))}{g(G^{-1}(u))}\nonumber\\
			=&\Lambda_Y^{w_2} (u).\nonumber
		\end{align}
		Now using (\ref{GWCPJdef}), 
		\begin{align*}
			\bar{\xi}^{w_1}(X)&= -\frac{1}{2}E\left(\Lambda_X^{w_1} (U)\right)\\
			&\ge -\frac{1}{2}E\left(\Lambda_Y^{w_2} (U)\right)\\
			&=\bar{\xi}^{w_2}(Y).
		\end{align*}
		(b) Proof is Similar to part (a). \hfill $\blacksquare$
		
		If we take $w_1(x)=w_2(x)=w(x)$ in the above theorem, then we have the following corollary.
		
		\begin{corollary}
			Let $X$  and $Y$ be nonnegative random variables with pdf's $f$ and $g$, cdf's $F$ and $G$, respectively having $u_X=u_Y<\infty$. Let $w$ is decreasing. \\
			(a) If  $X\le_{disp} Y$, then $\bar{\xi}^{w}(X)\ge \bar{\xi}^{w}(Y)$.\\
			(b) If  $X\ge_{disp} Y$, then $\bar{\xi}^{w}(X)\le \bar{\xi}^{w}(Y)$.
		\end{corollary}
		
		In the following theorem, we will study the behaviour of GWCRJ with respect to weight functions $w_1$ and $w_2$ and the dispersive order of random variables $X$ and $Y$.

		\begin{theorem}\label{thm com en1}
			Let $X$  and $Y$ be nonnegative random variables with pdf's $f$ and $g$, cdf's $F$ and $G$, respectively having $u_X=u_Y<\infty$.\\
			\\
			(a) If $w_1$ is decreasing, $w_1(x)\leq w_2(x)$ and $X\le_{disp} Y$, then ${\xi}^{w_1}(X)\ge {\xi}^{w_2}(Y)$.\\
			(b)   If $w_1$ is decreasing, $w_1(x)\geq w_2(x)$ and $X\ge_{disp} Y$, then ${\xi}^{w_1}(X)\le {\xi}^{w_2}(Y)$.
		\end{theorem}
		\noindent \textbf{Proof}
		(a) Since $X\le_{disp} Y$, therefore we have $f(F^{-1}(u))\ge g(G^{-1}(u))$ for $u \in (0,1)$. Then using Theorem 3.B.13(b) of Shaked and Shanthikumar (2007), $X\le_{disp} Y$ implies that $X\ge_{st} Y$. Hence $F^{-1}(u) \ge G^{-1}(u)$ for all $u\in (0,1)$. Since $w_1$ is decreasing and $w_1(x)\leq w_2(x)$ , then $w_1(F^{-1}(u)) \le w_1(G^{-1}(u))\le w_2(G^{-1}(u))$. 
		Hence 
		\begin{align}
			\Delta_X^{w_1} (u)&= \frac{ (1-u)^2 w_1(F^{-1}(u))}{f(F^{-1}(u))}\nonumber\\
			&\le  \frac{(1-u)^2 w_2(G^{-1}(u))}{g(G^{-1}(u))}\nonumber\\
			=&\Delta_Y^{w_2} (u).
		\end{align}
		Now using (\ref{GWCRJdef}), 
		\begin{align*}
			{\xi}^{w_1}(X)&= -\frac{1}{2}E\left(\Delta_X^{w_1} (U)\right)\\
			&\ge -\frac{1}{2}E\left(\Delta_Y^{w_2} (U)\right)\\
			&={\xi}^{w_2}(Y).
		\end{align*}
		(b) Proof is Similar to part (a). \hfill $\blacksquare$
		
		If we take $w_1(x)=w_2(x)=w(x)$ in the above theorem, then we have the following corollary.
		
		\begin{corollary}
			Let $X$  and $Y$ be nonnegative random variables with pdf's $f$ and $g$, cdf's $F$ and $G$, respectively having $u_X=u_Y<\infty$. Let $w$ is decreasing. \\
			(a) If  $X\le_{disp} Y$, then ${\xi}^{w}(X)\ge {\xi}^{w}(Y)$.\\
			(b) If  $X\ge_{disp} Y$, then ${\xi}^{w}(X)\le {\xi}^{w}(Y)$.
		\end{corollary}

		\section{Simple Random Sampling} \label{section3}

		Let $X$ be a random variable with finite mean $\mu$ and variance $\sigma^2$. For $\textbf{X}_{SRS}=\{X_i,\ i=1,\ldots,n\}$, the joint pdf is $\prod_{i=1}^{n}f(x_i)$, as $X_i$'s, $i=1,\ldots,n$ are independent and identically distributed (i.i.d.). Hence the GWCPJ  of $\textbf{X}_{SRS}^{(n)}$ can be defined as
		\begin{align}\label{srsxiw}
			\bar{\xi}^{w}(\textbf{X}_{SRS}^{(n)})&=\frac{-1}{2}\prod_{i=1}^{n}\left(\int_{-\infty}^{\infty}w(x_i)F^2(x_i)dx_i\right)\nonumber \\
			&=\frac{-1}{2}\left(-2\bar{\xi}^{w}(X)\right)^n\nonumber \\
			&=\frac{-1}{2}\left(\mathbb{E}(\Lambda_X^{w}(U))\right)^n.
		\end{align}
		
		\noindent Also, the GWCRJ  of $\textbf{X}_{SRS}^{(n)}$ can be defined as
		\begin{align}
			{\xi}^{w}(\textbf{X}_{SRS}^{(n)})&=\frac{-1}{2}\prod_{i=1}^{n}\left(\int_{-\infty}^{\infty}w(x_i)\bar{F}_X^2(x_i)dx_i\right)\nonumber \\
			&=\frac{-1}{2}\left(-2{\xi}^{w}(X)\right)^n\nonumber \\
			&=\frac{-1}{2}\left(\mathbb{E}(\Delta_X^{w}(U))\right)^n, \nonumber 
		\end{align}
		where $\Lambda_X^{w}(U)$ and $\Delta_X^{w}(U)$ are given in section \ref{section1}.
		\begin{example}
			If U is a uniform random variable on interval (0,1)  and $w(x)=x^m, m >0 $ then 
			\begin{align*}
				&\mathbb{E}[	\Lambda_X^{w}(U)] 		 
				= \mathbb{E}\left[ {U^2 w(U)}\right] = \frac{1}{m+3},\\
				&\mathbb{E}[ \Delta_X^{w}(U)] 
				= \mathbb{E}\left[ {(1-U)^2 w(U)}\right] = \frac{1}{(m+1)}-\frac{2}{(m+2)}+\frac{1}{(m+3)}, \\
				&\bar{\xi}^{w}(\textbf{X}_{SRS}^{(n)})
				=\frac{-1}{2}\left(\frac{1}{m+3}\right)^n, \\
				&{\xi}^{w}(\textbf{X}_{SRS}^{(n)})
				=\frac{-1}{2}\left(\frac{1}{(m+1)}-\frac{2}{(m+2)}+\frac{1}{(m+3)}\right)^n \nonumber.
			\end{align*} 		
			
		\end{example}
		
		\begin{theorem}\label{thm3sec3}
			Let $X$ be a non-negative absolutely continuous random variable with pdf $f$ and cdf $F.$ Assume $\psi(x)$ is an increasing function and  ${w(\psi(x))}{\psi^\prime (x)} \leq (\geq) w(x)$ and $\psi(0)=0$. If $Y=\psi(X)$, then $\bar{\xi}^{w}(\textbf{X}_{SRS}^{(n)}) \leq (\geq) \bar{\xi}^{w}(\textbf{Y}_{SRS}^{(n)})$.
		\end{theorem} 
		\noindent \textbf{Proof} Let the random variable $Y$ has pdf $g$ and cdf $G.$ Then 
		\begin{align*}
			\Lambda_Y^{w}(u)&=\frac{u^2 w\left(G^{-1}(u)\right)}{g\left(G^{-1}(u)\right)} =\frac{u^2 w(\psi(F^{-1}(u)))\psi^\prime(F^{-1}(u))}{f(F^{-1}(u))} \hspace{5mm} \forall \hspace{5mm} 0<u<1.
		\end{align*}
		Note that $\psi(x)\geq \psi(0)$, $\forall$ $x\geq 0$. Since ${w(\psi(x))}{\psi^\prime (x)} \leq (\geq) w(x).$ Hence for $ 0<u<1$, we have
		\begin{align*}
			\Lambda_Y^{w}(u)=\frac{u^2 w(\psi(F^{-1}(u)))\psi^\prime(F^{-1}(u))}{f(F^{-1}(u))} \leq (\geq) \frac{u^2{w((F^{-1}(u)))}}{f(F^{-1}(u))}=\Lambda_X^{w}(u).
		\end{align*}
		Therefore $\bar{\xi}^{w}(\textbf{X}_{SRS}^{(n)}) \geq(\leq) \bar{\xi}^{w}(\textbf{Y}_{SRS}^{(n)})$ using equation (\ref{srsxiw}). \hfill $\blacksquare$
		
		\begin{example}\label{eg2sec3}			
			Let $\psi(x)=e^x-1,\ \ w(x)=x, \  x\geq 0.$ It is easy to verify that $\psi(x)$ is an increasing function and $\psi(0)=0.$ Let $h(x)={w(\psi(x))}{\psi^\prime (x)}-w(x).$ 	That is, $h(x)= e^{2x}-e^x-x.$ Note that $h(0)=h^{'}(0)=0.$
			Moreover, $h^{'}(x)\geq 0 \ \forall \ x\geq 0.$ Thus, $h(x)$ is increasing function on $x \in [0,\infty].$ So, $h(x)\geq h(0)=0,$ that is, $w(\psi(x)){\psi^\prime (x)} \geq w(x) \ \forall x \geq 0.$	Therefore, using Theorem \ref{thm3sec3} for $Y=\psi(X)$, we have,
			\[\bar{\xi}^{w}(\textbf{X}_{SRS}^{(n)}) \geq \bar{\xi}^{w}(\textbf{Y}_{SRS}^{(n)}).\]
		\end{example}
		
		\section{MinRSSU and maxRSSU for GWCRJ and GWCPJ }\label{section4}
		Let $X_{i:i}$ represent $i$-th order statistics from a sample of size $i$ and $X_{1:i}$ represent $1$-st order statistics from a sample of size $i$. Let $F_{i:i}$ denote cdf of $i$-th order statistics from a sample of size $i$ and $F_{1:i}$ denotes cdf of $1$-th order statistics from a sample of size $i$. 
		
		The GWCPJ  of $\textbf{X}_{maxRSSU}^{(n)}$ can be defined as
		\begin{align}
			\bar{\xi}^{w}(\textbf{X}_{maxRSSU}^{(n)})
			&=\frac{-1}{2} \prod_{i=1}^{n} \left[-2{\bar \xi}^{w}(X_{i:i})\right]\nonumber \\
			&=\frac{-1}{2} \prod_{i=1}^{n} \left[ \int_{-\infty}^{\infty}w(x){F^2_{i:i}}(x)dx \right]\nonumber \\
			&=\frac{-1}{2} \prod_{i=1}^{n} \left[ \int_{-\infty}^{\infty}w(x){F^{2i}_X}(x)dx \right]\nonumber\\
			&=\frac{-1}{2} \prod_{i=1}^{n} \left[ \int_{0}^{1}\frac{{u}^{2i} w(F^{-1}(u))}{f(F^{-1}(u))}du \right]\nonumber\\
			&= \frac{-1}{2} \prod_{i=1}^{n} \mathbb{E} \left[ \frac{{U}^{2i} w(F^{-1}(U))}{f(F^{-1}(U))} \right] \nonumber  \\
			&= \frac{-1}{2} \prod_{i=1}^{n} \mathbb{E} \left[ \Psi^w_X(U) \right], \label{xibarmaxRSS},
		\end{align}
		where $\Psi^w_X(u)=\frac{{u}^{2i} w(F^{-1}(u))}{f(F^{-1}(u))}.$

		\noindent The GWCRJ  of $\textbf{X}_{minRSSU}^{(n)}$ can be defined as
		\begin{align}
			{\xi}^{w}(\textbf{X}_{minRSSU}^{(n)})
			&=\frac{-1}{2} \prod_{i=1}^{n} \left[-2{ \xi}^{w}(X_{1:i})\right]\nonumber \\
			&=\frac{-1}{2} \prod_{i=1}^{n} \left[ \int_{-\infty}^{\infty}w(x){\bar{F}_{1:i}}^2(x)dx \right]\nonumber \\
			&=\frac{-1}{2} \prod_{i=1}^{n} \left[ \int_{-\infty}^{\infty}w(x){\bar{F}}^{2i}(x)dx \right]\nonumber\\
			&=\frac{-1}{2} \prod_{i=1}^{n} \left[ \int_{0}^{1}\frac{{(1-u)}^{2i} w(F^{-1}(u))}{f(F^{-1}(u))}du \right]\nonumber\\
			&= \frac{-1}{2} \prod_{i=1}^{n} \mathbb{E} \left[ \frac{{(1-U)}^{2i} w(F^{-1}(U))}{f(F^{-1}(U))} \right] \nonumber \\
			&= \frac{-1}{2} \prod_{i=1}^{n} \mathbb{E} \left[ \Phi^w_X(U) \right], \label{ximinRSS},
		\end{align}
		where $\Phi^w_X(u)=\frac{{(1-u)}^{2i} w(F^{-1}(u))}{f(F^{-1}(u))}.$
		
		\noindent Let us see some examples to illustrate the above measures.
		\begin{example}
			If U is a uniform random variable on interval (0,1)  and $w(x)=x^m, m >0 $ then
			\begin{align*}
				&\bar{\xi}^{w}(\textbf{X}_{maxRSSU}^{(n)})
				= \frac{-1}{2} \prod_{i=1}^{n} \left(\frac{1}{2i+m+1}\right), \\
				&	{\xi}^{w}(\textbf{X}_{minRSSU}^{(n)})
				= -\frac{(\Gamma(m+1))^2}{2} \prod_{i=1}^{n} \frac{\Gamma(2i+1)}{\Gamma(2i+m+2)},
			\end{align*} 		
			
			where \[\Gamma(\alpha)=\int_{0}^{\infty} e^{-x} x^{\alpha-1} dx\].
		\end{example}
		
		\begin{example}
			Let $X$ be an exponential distribution with cdf $F(x)=1-e^{-\lambda x}, \ \lambda >0, \ x>0$. Let $w(x)=x^m, \ m>0, \ x>0$, then 
			\begin{align*}
				{\xi}^{w}(\textbf{X}_{minRSSU}^{(n)})
				&=\frac{-1}{2} \prod_{i=1}^{n} \left[ \int_{0}^{1}\frac{{(1-u)}^{2i} w(F^{-1}(u))}{f(F^{-1}(u))}du \right]\\
				&= \frac{-1}{2} \prod_{i=1}^{n} \frac{(-1)^m}{\lambda^{m+1}} \int_{0}^{1} (1-u)^{2i-1} (\log(1-u))^m du. 
			\end{align*} 
			Taking $1-u=e^{-x}$ in above expression, we get 
			
			\begin{align*}
				{\xi}^{w}(\textbf{X}_{minRSSU}^{(n)})
				&= -\frac{1}{2} \prod_{i=1}^{n} \frac{1}{\lambda^{m+1}} \int_{0}^{\infty} x^m e^{-2ix}  dx\\
				&= -\frac{1}{2} \prod_{i=1}^{n} \frac{1}{\lambda^{m+1}} \frac{\Gamma(m+1)}{(2i)^{m+1}} dx\\
				&= -\frac{1}{2} \left(\frac{\Gamma(m+1)}{\lambda^{m+1}}\right)^n \prod_{i=1}^{n} \left(\frac{1}{2i}\right)^{m+1} \\
				&=-\frac{1}{2} \left(\frac{\Gamma(m+1)}{(2\lambda)^{m+1}}\right)^n \left(\frac{1}{n!}\right)^{m+1}.	
			\end{align*} 
		\end{example}
		\begin{example}
			Let a random variable $X$ have survival function $\bar{F}(x)=(1-x)^b, \ \ 0<x<1,  \ \ b>0. $ If $w(x)=x^{m}, \ \ x>0, \ \ m>0$ then 
			\begin{align*}
				{\xi}^{w}(\textbf{X}_{minRSSU}^{(n)})&= \frac{-1}{2} \prod_{i=1}^{n} \left[ \int_{0}^{1}w(x){\bar{F}}^{2i}(x)dx \right] \\
				&= \frac{-1}{2} \prod_{i=1}^{n} \left[ \int_{0}^{1} x^m {(1-x)}^{2ib}dx \right]\\
				&= \frac{-1}{2} \prod_{i=1}^{n} B(m+1, 2ib+1),     			
			\end{align*}
			where $B(\alpha, \beta)$ is beta function defined as  \[ B(\alpha, \beta)= \int_{0}^{1} x^{\alpha-1} (1-x)^{1-\beta}dx. \]
		\end{example}

		\begin{theorem}\label{new1}
			Let $X$ be a non-negative absolutely continuous random variable with pdf $f$ and cdf $F.$ Assume $\psi(x)$ is an increasing function and  ${w(\psi(x))}{\psi^\prime (x)} \leq (\geq) w(x)$ and $\psi(0)=0$. If \ $Y=\psi(X)$, then $	\bar{\xi}^{w}(\textbf{X}_{maxRSSU}^{(n)})\leq (\geq) 	\bar{\xi}^{w}(\textbf{Y}_{maxRSSU}^{(n)}).$
		\end{theorem}
		\noindent \textbf{Proof} Let random variable $Y$  has pdf $g$ and cdf $G.$ Since ${w(\psi(x))}{\psi^\prime (x)} \leq (\geq)  w(x)$ and $\psi(0)=0$, therefore, 
		\begin{align*}
			\Psi^w_Y(u)&=\frac{u^{2i} w\left(G^{-1}(u)\right)}{g\left(G^{-1}(u)\right)}\\
			&=\frac{u^{2i} w(\psi(F^{-1}(u)))\psi^\prime(F^{-1}(u))}{f(F^{-1}(u))} \\
			&\leq (\geq)  \frac{{u}^{2i} w(F^{-1}(u))}{f(F^{-1}(u))}=\Psi^w_X(u) .
		\end{align*}
		Hence, From (\ref{xibarmaxRSS}) we obtain
		\begin{align*}
			\bar{\xi}^{w}(\textbf{X}_{maxRSSU}^{(n)})\geq (\leq) 	\bar{\xi}^{w}(\textbf{Y}_{maxRSSU}^{(n)}). 
	\end{align*}}
	\hfill $\blacksquare$
	
	\begin{example}
		Let $\psi(x)=e^x-1,\ \ w(x)=x, \  x  \geq 0.$ Note that from example \ref{eg2sec3}, $\psi(x)$ is an increasing function and  ${w(\psi(x))}{\psi^\prime (x)} \geq w(x)$ and $\psi(0)=0.$ Using Theorem \ref{thm3sec3} for $Y=\psi(X),$ we have,
		\[ \bar{\xi}^{w}(\textbf{X}_{maxRSSU}^{(n)}) \geq 	\bar{\xi}^{w}(\textbf{Y}_{maxRSSU}^{(n)}) \].
	\end{example}
	
	\begin{theorem}
		Let $X$ be a non-negative absolutely continuous random variable with pdf $f$ and cdf $F.$ Assume $\psi(x)$ is an increasing function and  ${w(\psi(x))}{\psi^\prime (x)} \leq (\geq) w(x)$ and $\psi(0)=0$. If $Y=\psi(X)$, then $	{\xi}^{w}(\textbf{X}_{minRSSU}^{(n)})\leq (\geq) 	{\xi}^{w}(\textbf{Y}_{minRSSU}^{(n)}).$
	\end{theorem}
	
	\noindent \textbf{Proof} Let random variable $Y$  has pdf $g$ and cdf $G.$ Since ${w(\psi(x))}{\psi^\prime (x)} \leq (\geq)  w(x)$ and $\psi(0)=0$, therefore, 
	\begin{align*}
		\Phi^w_Y(u)&=\frac{(1-u)^{2i} w\left(G^{-1}(u)\right)}{g\left(G^{-1}(u)\right)}\\
		&=\frac{(1-u)^{2i} w(\psi(F^{-1}(u)))\psi^\prime(F^{-1}(u))}{f(F^{-1}(u))} \\
		&\leq (\geq)  \frac{{(1-u)}^{2i} w(F^{-1}(u))}{f(F^{-1}(u))}=\Phi^w_X(u) .
	\end{align*}
	Hence, From (\ref{ximinRSS}) we obtain
	\begin{align*}
		{\xi}^{w}(\textbf{X}_{minRSSU}^{(n)})\geq (\leq) 	{\xi}^{w}(\textbf{Y}_{minRSSU}^{(n)}). 
	\end{align*}
	\hfill $\blacksquare$
	
	\begin{example}
		Let $\psi(x)=e^x-1,\ \ w(x)=x, \  x  \geq 0.$ Note that from Example \ref{eg2sec3}, $\psi(x)$ is an increasing function and  ${w(\psi(x))}{\psi^\prime (x)} \geq w(x)$ and $\psi(0)=0.$ Using Theorem \ref{thm3sec3} for $Y=\psi(X),$ we have,
		\[  {\xi}^{w}(\textbf{X}_{minRSSU}^{(n)}) \geq 	{\xi}^{w}(\textbf{Y}_{minRSSU}^{(n)}) \].
	\end{example}
	
	\begin{theorem}
		Let $X$ be an absolutely continuous random variable with pdf f and cdf F. Then for $n\geq 2$,
	\end{theorem}
	\begin{align*}
		\bar{\xi}^{w}(\textbf{X}_{maxRSSU}^{(n)}) \geq \bar{\xi}^{w}(\textbf{X}_{SRS}^{(n)}).
	\end{align*}
	\noindent \textbf{Proof} Consider
	\begin{align*}
		\bar{\xi}^{w}(\textbf{X}_{maxRSSU}^{(n)}) &=\frac{-1}{2} \prod_{i=1}^{n} \mathbb{E} \left[ \frac{{U}^{2i} w(F^{-1}(U))}{f(F^{-1}(U))} \right] \\
		&=\frac{-1}{2} \prod_{i=1}^{n} \left[ \int_{0}^{1}\frac{{u}^{2i} w(F^{-1}(u))}{f(F^{-1}(u))}du \right] \\
		& \geq \frac{-1}{2} \prod_{i=1}^{n} \left[ \int_{0}^{1}\frac{{u}^{2} w(F^{-1}(u))}{f(F^{-1}(u))}du \right] \\
		&= \frac{-1}{2} \prod_{i=1}^{n} \mathbb{E} \left[ \frac{{U}^{2} w(F^{-1}(U))}{f(F^{-1}(U))} \right] \\
		&= \bar{\xi}^{w}(\textbf{X}_{SRS}^{(n)}).
	\end{align*}
	Hence the result.		\hfill $\blacksquare$
	
	\begin{theorem}
		Let $X$ be an absolutely continuous random variable with pdf f and cdf F. Then for $n\geq 2$,
	\end{theorem}
	\begin{align*}
		{\xi}^{w}(\textbf{X}_{minRSSU}^{(n)}) \geq {\xi}^{w}(\textbf{X}_{SRS}^{(n)}).
	\end{align*}
	\noindent \textbf{Proof}  Consider
	\begin{align*}
		{\xi}^{w}(\textbf{X}_{minRSSU}^{(n)}) &=\frac{-1}{2} \prod_{i=1}^{n} \mathbb{E} \left[ \frac{{(1-U)}^{2i} w(F^{-1}(U))}{f(F^{-1}(U))} \right] \\
		&=\frac{-1}{2} \prod_{i=1}^{n} \left[ \int_{0}^{1}\frac{{(1-u)}^{2i} w(F^{-1}(u))}{f(F^{-1}(u))}du \right] \\
		& \geq \frac{-1}{2} \prod_{i=1}^{n} \left[ \int_{0}^{1}\frac{{(1-u)}^{2} w(F^{-1}(u))}{f(F^{-1}(u))}du \right] \\
		&= \frac{-1}{2} \prod_{i=1}^{n} \mathbb{E} \left[ \frac{{(1-U)}^{2} w(F^{-1}(U))}{f(F^{-1}(U))} \right] \\
		&= {\xi}^{w}(\textbf{X}_{SRS}^{(n)}).
	\end{align*}
	Hence the result.		\hfill $\blacksquare$

	\section{Stochastic Comparision} \label{section5}
	
	The following lemma will be used in deriving Theorem \ref{theorem1usinglemma1} and Theorem \ref{theorem2usinglemma1}. Lemma 1 is also used by Qiu and Raqab (2022) and Gupta and Chaudhary (2022) to derive their results.
	\begin{lemma}\label{lemma1} [Ahmed et al. (1986)]
		Let $X$ and $Y$  be nonnegative random variables with pdf's $f$ and $g$, respectively, satisfying $f(0)\ge g(0)>0$. If $X\le_{su}Y$ (or $X\le_{*}Y$ or $X\le_{c}Y$), then $X\le_{disp}Y$.
	\end{lemma}

	In the following theorem, we will study the behaviour maxRSSU of GWCPJ with respect to weight functions $w_1$ and $w_2$ and dispersive order of random variables $X$ and $Y$.
	\begin{theorem}\label{thm1sec5}
		Let $X$  and $Y$ be nonnegative random variables with pdf's $f$ and $g$, cdf's $F$ and $G$, respectively having $u_X=u_Y<\infty$.\\
		\\
		(a) If $w_1$ is decreasing, $w_1(x)\leq w_2(x)$ and $X\le_{disp} Y$, then $\bar{\xi}^{w_1}(\textbf{X}_{maxRSSU}^{(n)})\ge \bar{\xi}^{w_2}(\textbf{Y}_{maxRSSU}^{(n)})$.\\
		(b)   If $w_1$ is decreasing, $w_1(x)\geq w_2(x)$ and $X\ge_{disp} Y$, then $\bar{\xi}^{w_1}(\textbf{X}_{maxRSSU}^{(n)})\le \bar{\xi}^{w_2}(\textbf{Y}_{maxRSSU}^{(n)})$.
	\end{theorem}
	\noindent \textbf{Proof} (a) Since $X\le_{disp} Y$, therefore we have $f(F^{-1}(u))\ge g(G^{-1}(u))$ for $u \in (0,1)$. Then using Theorem 3.B.13(b) of Shaked and Shanthikumar (2007), $X\le_{disp} Y$ implies that $X\ge_{st} Y$. Hence, $F^{-1}(u) \ge G^{-1}(u)$ for all $u\in (0,1)$. Since $w_1$ is decreasing and $w_1(x)\leq w_2(x)$ , then $w_1(F^{-1}(u)) \le w_1(G^{-1}(u))\le w_2(G^{-1}(u))$. 
	Hence,
	\begin{align}
		\Psi^w_X(u)=	u^{2i} \frac{ w_1(F^{-1}(u))}{f(F^{-1}(u))} \le u^{2i} \frac{ w_2(G^{-1}(u))}{g(G^{-1}(u))}=\Psi^w_Y(u)\nonumber
	\end{align}
	Now using (\ref{xibarmaxRSS}),
	\begin{align*}
		\bar{\xi}^{w_1}(\textbf{X}_{maxRSSU}^{(n)})\ge \bar{\xi}^{w_2}(\textbf{Y}_{maxRSSU}^{(n)}).
	\end{align*}
	(b) Proof is Similar to part (a). \hfill $\blacksquare$

	If we take $w_1(x)=w_2(x)=w(x)$ in the above theorem, then we have the following corollary.
	
	\begin{corollary}\label{cor1sec5}
		Let $X$  and $Y$ be nonnegative random variables with pdf's $f$ and $g$, cdf's $F$ and $G$, respectively having $u_X=u_Y<\infty$. Let $w$ is decreasing. \\
		(a) If  $X\le_{disp} Y$, then $\bar{\xi}^{w}(\textbf{X}_{maxRSSU}^{(n)})\ge \bar{\xi}^{w}(\textbf{Y}_{maxRSSU}^{(n)})$.\\
		(b) If  $X\ge_{disp} Y$, then $\bar{\xi}^{w}(\textbf{X}_{maxRSSU}^{(n)})\le \bar{\xi}^{w}(\textbf{Y}_{maxRSSU}^{(n)})$.
	\end{corollary}

	One may refer Shaked and Shantikumar (2007) for details of convex transform order ($\leq_c$), star order ($\leq_{\star}$), super additive order ($\leq_{su}$), and dispersive order ($\leq_{disp}$). In view of Theorem \ref{thm1sec5}  and Lemma \ref{lemma1}, the following result is obtained.

	\begin{theorem}\label{theorem1usinglemma1}
		Let $X$  and $Y$ be nonnegative random variables with pdf's $f$ and $g$, cdf's $F$ and $G$, respectively having $u_X=u_Y<\infty$.\\
		\\
		(a) If $w_1$ is decreasing, $w_1(x)\leq w_2(x)$ and $X\le_{su} Y$( or $X\le_{*} Y$ or $X\le_{c} Y$), then $\bar{\xi}^{w_1}(\textbf{X}_{maxRSSU}^{(n)})\ge \bar{\xi}^{w_2}(\textbf{Y}_{maxRSSU}^{(n)})$.\\
		(b)   If $w_1$ is decreasing, $w_1(x)\geq w_2(x)$ and $X\ge_{su} Y$ (or $X\ge_{*} Y$ or $X\ge_{c} Y$ ), then $\bar{\xi}^{w_1}(\textbf{X}_{maxRSSU}^{(n)})\le \bar{\xi}^{w_2}(\textbf{Y}_{maxRSSU}^{(n)})$.
	\end{theorem}
	
	If we take $w_1(x)=w_2(x)=w(x)$ in the above theorem, then we have the following corollary.
	
	\begin{corollary}
		Let $X$  and $Y$ be nonnegative random variables with pdf's $f$ and $g$, cdf's $F$ and $G$, respectively having $u_X=u_Y<\infty$.\\
		\\
		(a) If $w$ is decreasing and $X\le_{su} Y$( or $X\le_{*} Y$ or $X\le_{c} Y$), then $\bar{\xi}^{w}(\textbf{X}_{maxRSSU}^{(n)})\ge \bar{\xi}^{w}(\textbf{Y}_{maxRSSU}^{(n)})$.\\
		(b)   If $w$ is decreasing and $X\ge_{su} Y$ (or $X\ge_{*} Y$ or $X\ge_{c} Y$ ), then $\bar{\xi}^{w}(\textbf{X}_{maxRSSU}^{(n)})\le \bar{\xi}^{w}(\textbf{Y}_{maxRSSU}^{(n)})$.
	\end{corollary}

	In the following theorem, we will study the behaviour of minRSSU of GWCRJ with respect to weight functions $w_1$ and $w_2$ and the dispersive order of random variables $X$ and $Y$.

	\begin{theorem}\label{thm2sec5}
		Let $X$  and $Y$ be nonnegative random variables with pdf's $f$ and $g$, cdf's $F$ and $G$, respectively having $u_X=u_Y<\infty$.\\
		\\
		(a) If $w_1$ is decreasing, $w_1(x)\leq w_2(x)$ and $X\le_{disp} Y$, then ${\xi}^{w_1}(\textbf{X}_{minRSSU}^{(n)})\ge {\xi}^{w_2}(\textbf{Y}_{minRSSU}^{(n)})$.\\
		(b)   If $w_1$ is decreasing, $w_1(x)\geq w_2(x)$ and $X\ge_{disp} Y$, then ${\xi}^{w_1}(\textbf{X}_{minRSSU}^{(n)})\le {\xi}^{w_2}(\textbf{Y}_{minRSSU}^{(n)})$.
	\end{theorem}
	\noindent \textbf{Proof} (a) Since $X\le_{disp} Y$, therefore we have $f(F^{-1}(u))\ge g(G^{-1}(u))$ for $u \in (0,1)$. Then using Theorem 3.B.13(b) of Shaked and Shanthikumar (2007), $X\le_{disp} Y$ implies that $X\ge_{st} Y$. Hence $F^{-1}(u) \ge G^{-1}(u)$ for all $u\in (0,1)$. Since $w_1$ is decreasing and $w_1(x)\leq w_2(x)$ , then $w_1(F^{-1}(u)) \le w_1(G^{-1}(u))\le w_2(G^{-1}(u))$. 
	Hence
	\begin{align}
		\Phi^w_X(u)=(1-u)^{2i} \frac{ w_1(F^{-1}(u))}{f(F^{-1}(u))} \le (1-u)^{2i} \frac{ w_2(G^{-1}(u))}{g(G^{-1}(u))}=\Phi^w_Y(u)\nonumber
	\end{align}
	Now using (\ref{ximinRSS}),
	\begin{align*}
		{\xi}^{w_1}(\textbf{X}_{minRSSU}^{(n)})\ge {\xi}^{w_2}(\textbf{Y}_{minRSSU}^{(n)}).
	\end{align*}
	(b) Proof is Similar to part (a). \hfill $\blacksquare$

	If we take $w_1(x)=w_2(x)=w(x)$ in the above theorem, then we have the following corollary.
	
	\begin{corollary}
		Let $X$  and $Y$ be nonnegative random variables with pdf's $f$ and $g$, cdf's $F$ and $G$, respectively having $u_X=u_Y<\infty$. Let $w$ is decreasing. \\
		(a) If  $X\le_{disp} Y$, then ${\xi}^{w}(\textbf{X}_{minRSSU}^{(n)})\ge {\xi}^{w}(\textbf{Y}_{minRSSU}^{(n)})$.\\
		(b) If  $X\ge_{disp} Y$, then ${\xi}^{w}(\textbf{X}_{minRSSU}^{(n)})\le {\xi}^{w}(\textbf{Y}_{minRSSU}^{(n)})$.
	\end{corollary}
	
	In view of Theorem \ref{thm2sec5}  and Lemma \ref{lemma1}, the following result is obtained.

	\begin{theorem}\label{theorem2usinglemma1}
		Let $X$  and $Y$ be nonnegative random variables with pdf's $f$ and $g$, cdf's $F$ and $G$, respectively having $u_X=u_Y<\infty$.\\
		\\
		(a) If $w_1$ is decreasing, $w_1(x)\leq w_2(x)$ and $X\le_{su} Y$ (or $X\le_{*} Y$ or $X\le_{c} Y$), then ${\xi}^{w_1}(\textbf{X}_{minRSSU}^{(n)})\ge {\xi}^{w_2}(\textbf{Y}_{minRSSU}^{(n)})$.\\
		(b)   If $w_1$ is decreasing, $w_1(x)\geq w_2(x)$ and $X\ge_{su} Y$ (or $X\ge_{*} Y$ or $X\ge_{c} Y$), then ${\xi}^{w_1}(\textbf{X}_{minRSSU}^{(n)})\le {\xi}^{w_2}(\textbf{Y}_{minRSSU}^{(n)})$.
	\end{theorem}
	
	If we take $w_1(x)=w_2(x)=w(x)$ in the above theorem, then we have the following corollary.

	\begin{corollary}
		Let $X$  and $Y$ be nonnegative random variables with pdf's $f$ and $g$, cdf's $F$ and $G$, respectively having $u_X=u_Y<\infty$.\\
		\\
		(a) If $w$ is decreasing and $X\le_{su} Y$ (or $X\le_{*} Y$ or $X\le_{c} Y$), then ${\xi}^{w}(\textbf{X}_{minRSSU}^{(n)})\ge {\xi}^{w}(\textbf{Y}_{minRSSU}^{(n)})$.\\
		(b)   If $w$ is decreasing and $X\ge_{su} Y$ (or $X\ge_{*} Y$ or $X\ge_{c} Y$), then ${\xi}^{w}(\textbf{X}_{minRSSU}^{(n)})\le {\xi}^{w}(\textbf{Y}_{minRSSU}^{(n)})$.
	\end{corollary}

	\section{Monotone Properties } \label{section6}
	
	The following result gives the conditions under which the ${\xi}^{w}(\textbf{X}_{minRSSU}^{(n)})$  will increase with respect to $n$.

	\begin{theorem}
		If $\frac{w(F^{-1}(u)}{f(F^{-1}(u))}\leq  1, \ 0<u<1,\ $ then ${\xi}^{w}(\textbf{X}_{minRSSU}^{(n)})$ is increasing in $n\geq 1.$
	\end{theorem}
	\noindent \textbf{Proof} Consider
	\begin{align*} \frac{{\xi}^{w}(\textbf{X}_{minRSSU}^{(n+1)})}{{\xi}^{w}(\textbf{X}_{minRSSU}^{(n)})} =\int_{0}^{1}  \frac{(1-u)^{2n+2} w(F^{-1}(u))}{f(F^{-1}(u))}du \leq \frac{1}{2n+3}\leq 1 \ \forall \ n\geq1.
	\end{align*}
	Hence, \[{\xi}^{w}(\textbf{X}_{minRSSU}^{(n+1)}) \geq {\xi}^{w}(\textbf{X}_{minRSSU}^{(n)}). \]

	The following result gives the conditions under which the $\bar{\xi}^{w}(\textbf{X}_{maxRSSU}^{(n)})$  will increase with respect to $n$.

	\begin{theorem}
		If $\frac{w(F^{-1}(u)}{f(F^{-1}(u))}\leq 1, \ 0<u<1,\ $ then $\bar{\xi}^{w}(\textbf{X}_{maxRSSU}^{(n)})$ is increasing in $n\geq 1.$
	\end{theorem}
	\noindent \textbf{Proof} Consider
	\begin{align*} \frac{\bar{\xi}^{w}(\textbf{X}_{maxRSSU}^{(n+1)})}{\bar{\xi}^{w}(\textbf{X}_{maxRSSU}^{(n)})} = \int_{0}^{1}  \frac{u^{2n+2} w(F^{-1}(u))}{f(F^{-1}(u))}du \leq \frac{1}{2n+3}\leq 1 \ \forall \ n\geq1.
	\end{align*}
	Hence, \[\bar{\xi}^{w}(\textbf{X}_{maxRSSU}^{(n+1)}) \geq \bar{\xi}^{w}(\textbf{X}_{maxRSSU}^{(n)}). \]

	\section{Empirical measures of GWCPJ and GWCRJ}\label{section7}
	The empirical measure of $F$ is defined as 		
	\begin{eqnarray*}
		\reallywidehat{F}_{n}(x)=
		\begin{cases}
			0, \hspace{4mm} x<X_{1:n}\\
			\frac{i}{n}, \hspace{4mm} X_{i:n}\leq x<X_{i+1:n}, \ \ \ i=1,2,\dots, n-1.\\
			1, \hspace{5mm} x \geq X_{n:n}
		\end{cases}
	\end{eqnarray*}
	where $X_{1:n}\leq X_{2:n} \leq X_{3:n} \leq \dots \leq X_{n:n}$ denote order statistics of random sample $X_1, X_2, \dots, X_n.$

	The empirical measure of GWCPJ for $w(x)=x^m, \ x>0, \ m>0$ is
	\begin{align*}
		\reallywidehat {\bar{\xi}^m}(F_n)&=-\frac{1}{2} \int_{-\infty}^{\infty} x^m \reallywidehat{F^2_n}(x)dx = -\frac{1}{2} \sum_{i=1}^{n-1} \int_{x_{i:n}}^{x_{i+1:n}} x^m  \left(\frac{i}{n}\right)^2dx \\
		&=  -\frac{1}{2(m+1)}  \sum_{i=1}^{n-1}  \left[ (x^{m+1}_{i+1:n}-x^{m+1}_{i:n}) \left(\frac{i}{n}\right)^2 \right].
	\end{align*}

	The empirical measure of GWCRJ for $w(x)=x^m, \ x>0, \ m>0$ is
	\begin{align*}
		\reallywidehat {{\xi}^m}(F_n)&=-\frac{1}{2} \int_{-\infty}^{\infty} x^m \reallywidehat{\bar{F}^2_n}(x)dx = -\frac{1}{2} \sum_{i=1}^{n-1} \int_{x_{i:n}}^{x_{i+1:n}} x^m  \left(\frac{i}{n}\right)^2dx \\
		&=  -\frac{1}{2(m+1)}  \sum_{i=1}^{n-1}  \left[ (x^{m+1}_{i+1:n}-x^{m+1}_{i:n}) \left(1-\frac{i}{n}\right)^2 \right].
	\end{align*}
	
	Analogous to Theorem 9 of Rao et al. (2004), Theorem 5.1 of Tahmasebi et al. (2020) and Theorem 4.1 of Tahmasebi and Toomaj (2020), we state the following theorems and the proof follows similarly.

	\begin{theorem}
		Let $X$ be a non negative absolutely continuous random variable with cdf $F.$ Then for any random $X$ in $L^p$ for some $p>2,$     	$\reallywidehat {\bar{\xi}^m}(F_n)$ converges almost surely to ${\bar{\xi}^m}(X)$ as $n \  \rightarrow \infty.$
	\end{theorem}
	
	\begin{theorem}
		Let $X$ be a non negative absolutely continuous random variable with cdf $F.$ Then for any random $X$ in $L^p$ for some $p>2,$ 
		$\reallywidehat {{\xi}^m}(F_n)$ converges almost surely to ${{\xi}^m}(X)$ as $n \  \rightarrow \infty.$
	\end{theorem}
	
	We create another estimator based on the kernel-smoothed estimator of the distribution function because smoothed estimators perform better than non-smoothed estimators.

	The second estimator of GWCPJ for $w(x)=x^m, \ x>0, \ m>0$ is
	\begin{align*}
		\reallywidehat {\bar{\xi}^m}(F_h)&=-\frac{1}{2} \int_{-\infty}^{\infty} x^m {F^2_h}(x)dx = -\frac{1}{2} \sum_{i=1}^{n-1} \int_{x_{i:n}}^{x_{i+1:n}} x^m  \left({F^2_h}(x_i) \right)^2dx \\
		&=  -\frac{1}{2(m+1)}  \sum_{i=1}^{n-1}  \left[ (x^{m+1}_{i+1:n}-x^{m+1}_{i:n}) \left( {F^2_h}(x_i)\right)^2 \right].
	\end{align*}
	The second estimator of GWCRJ for $w(x)=x^m, \ x>0, \ m>0$ is
	\begin{align*}
		\reallywidehat {{\xi}^m}(F_h)&=-\frac{1}{2} \int_{-\infty}^{\infty} x^m {\bar{F}^2_h}(x)dx = -\frac{1}{2} \sum_{i=1}^{n-1} \int_{x_{i:n}}^{x_{i+1:n}} x^m  \left({\bar{F}^2_h}(x_i) \right)dx \\
		&=  -\frac{1}{2(m+1)}  \sum_{i=1}^{n-1}  \left[ (x^{m+1}_{i+1:n}-x^{m+1}_{i:n}) \left({\bar{F}^2_h}(x_i)\right) \right],
	\end{align*}
	where $\bar{F}=1-F$ and $F_h$ is defined by Nadaraya (1964) as 
	\[F_h(x)=\frac{1}{n} \sum_{i=1}^{n} L \left(\frac{x-X_i}{h}\right).\]
	Here, $L(x)=\int_{-\infty}^{x} K(t) dt$, that is, $L$ is cdf of positive kernel $K$ and $h$ is a bandwidth parameter (see Sarda (1993) and Jahanshahi (2020)).

	\section{Conclusion} \label{s10conclusion}\label{section8}
	Generalized versions of cumulative past and residual extropy have been defined. The general GWCRJ and GWCPJ uncertainty measures based on the minRSSU and maxRSSU data, respectively, have been introduced in this study. Additionally, this work presented the GWCRJ and GWCPJ-based uncertainty measures of SRS. For minRSSU, maxRSSU, and SRS data, some results of the GWCPJ and GWCRJ measures were obtained, including stochastic orders. With help of empirical measure of the distribution function, we obtained two estimators of GWCRJ and GWCPJ.\\

	\noindent \textbf{\Large Funding} \\
	\\
	Santosh Kumar Chaudhary would like to thank the council of scientific and industrial research (CSIR), Government of India (File Number 09/0081 (14002)/2022-EMR-I) for financial assistance.\\
	\\
	\textbf{ \Large Conflict of interest} \\
	\\
	The authors declare no conflict of interest.

	\vspace{.3in}
	\noindent
	{\bf Santosh Kumar Chaudhary}\\
	Department of Mathematics,\\
	Indian Institute of Technology Kharagpur\\
	Kharagpur-721302, INDIA\\
	E-mail: skchaudhary1994@kgpian.iitkgp.ac.in

	\vspace{.1in}
	
	\noindent
	{\bf  Nitin Gupta} \\
	Department of Mathematics,\\
	Indian Institute of Technology Kharagpur\\
	Kharagpur-721302, INDIA\\
	E-mail: nitin.gupta@maths.iitkgp.ac.in	\\
\end{document}